\documentclass[titlepage,twoside,12pt]{article}
\usepackage{amssymb}
\usepackage{amsfonts}
\begin{document}

{\Large \bf The Nash-Equilibrium Requires \\ \\ Strong Cooperation} \\ \\

{\it Elem\'{e}r E Rosinger \\ Department of Mathematics \\ University of Pretoria \\ Pretoria, 0002 South Africa \\
e-mail : eerosinger@hotmail.com} \\ \\

{\bf Abstract} \\

Contrary to the customary view that the celebrated Nash-equilibrium theorem in Game Theory is
paradigmatic for non-cooperative games, it is shown that, in fact, it is essentially based on
a particularly strong cooperation assumption. Furthermore, in practice, this cooperation
assumption is simply unrealistic. \\ \\

{\large \bf 1. Introduction} \\

One of the major divides in Game Theory is between {\it cooperative} and {\it
non-cooperative} games. And as it happened along the development of that theory around the
times of World War II, the idea emerged that, in certain important ways, non-cooperative games
are more basic than the cooperative ones. In a few words, the thinking leading to that idea
could be formulated as follows. A cooperative game between $n \geq 2$ players, if it is to be
nontrivial, must leave each player some freedom, no matter which may be the binding agreements
of the respective cooperation. And then, if one focuses on the respective freedom of each of
the players, one may obtain a game which is like a non-cooperative one. \\

The difficulty with such a argument is that games, be they cooperative or non-cooperative, can
be utterly complex. And to mention but one such aspect of complexity, let us recall that many
games prove to be {\it algorithmically unsolvable}, Binmore [1-3]. \\

A consequence of the significant complexities involved in games is that it may not always be
so easy to separate or subtract the cooperative part of a game, in order to remain with a game
which may be considered as being non-cooperative. In other words, the formula \\

$~~~~~~ game ~~\equiv~~ non-cooperative~game~~~ ( modulo~~ cooperation ) $ \\

does {\it not} in general hold. Therefore, one is not justified in seeing non-cooperative
games as being more basic than the cooperative ones.  \\

Regardless of the above, one may note that the centrality in Game Theory of the
Nash-equilibrium theorem about non-cooperative games is due not in a small measure to a
certain tacit or lingering acceptance of the above incorrect formula. \\

The aim of this paper, however, is not so much about the clarification of the place of the
Nash-equilibrium theorem with Game Theory, but rather to show that, no matter how much that
theorem is seen as a {\it non-cooperative paradigm}, it does essentially depend on an
assumption which requires such a {\it strong cooperation}, as to render it {\it unrealistic}
from practical point of view. \\

Before going further, it may be useful to recall in short certain formative moments in the
evolution of modern Game Theory, moments which may shed a light upon the mentioned views
relating to the divide between cooperative and non-cooperative games. \\

Game theory, as initiated by John von Neumann, see von Neumann \& Morgenstern, centers around
the {\it individual} and aims to lead to a rational outcome when two or more individuals
interact in well defined situations and do so, however, without the interference of an overall
authority. The very development since the 1920s of the mathematical theory of games is in
itself an act of rational behaviour, albeit on a certain {\it meta-level}, which involves both
the levels of the interest in developing the general concepts, axioms, theorems, and so on, as
well as the levels at which they are put to their effective uses in a variety of relevant
applications. \\
And yet, such a {\it meta-rational} behaviour appears to have mostly come to a halt when
dealing with {\it cooperation}. \\
Instead, the effort has rather been focused upon non-cooperative contexts, and rationality got
thus limited to them. A further aggravation of such a limitation upon rationality has come
from the fact that cooperation is so often being associated with issues of ethics, morality,
wisdom, philosophy, or on the contrary, politics, or mere expediency. \\

One of the few more prominent contributions to cooperation has been the 1984 study of Robert
Axelrod, which however is limited to two person non zero sum games, and it centers around one
of the simplest nontrivial such examples, namely, the celebrated game called the Prisoner's
Dilemma, suggested in the early 1950s by Merrill Flood and Melvin Drescher, and formalized by
Albert W Tucker, see Axelrod, Rasmusen. \\

However, even such simplest nontrivial examples can show that cooperation itself is often but
a most natural matter of rationality, albeit manifested in more evolved forms. Indeed, in many
non-cooperative situations what can be obtained by those involved proves to be {\it
significantly less} than what may be available through suitable cooperation. Thus the choice
of cooperation need not at all be seen as merely a matter which has to do with expediency,
politics, wisdom, morality or ethics. \\

And then the issue is simply the following : do we limit rationality and stop it before
considering cooperation in ways more adequate to its considerable depth and potential, or
instead, are we ready to try to be rational all the way ? \\

And the fact is that very few situations are of a nature in which competition or conflict is
total, and total to the extent that there simply cannot be any place whatsoever for one or
another form of cooperation. \\
In game theory such a situation corresponds to extreme cases such as the {\it two person zero
sum game}. \\
In rest, that is, in the vast majority of cases, the possibilities for cooperation and its
considerable rewards may often be there. And then, all it takes is to {\it extend} and {\it
deepen} rationality, and thus find and develop suitable ways of cooperation. \\

Needless to say, cooperative games are by no means less complex than the non-cooperative ones.
Therefore, one can expect yet another reason why attention may focus more on the latter, and
why one may try to reduce the former to them. \\
However, as follows from Binmore [1-3], this limitation of rationality mostly to the study of
noncooperative games is still leading to complexities which are not algorithmically
solvable. \\
During the late 1940s and early 1950s when game theory appears to have known its period of
most massive development, there was not much awareness about the possibility of the presence
of the type of deep difficulties which would more than three decades later be pointed out by
Binmore. \\

In such a context, during the late 1940s and early 1950s, the so called "Nash Program" emerged
which aimed to reduce cooperative games to non-cooperative ones, see Nash [1,3]. \\

There is, however, a rather weighty reason, explanation, and maybe also excuse, for the fact
that the Nash Program has known a certain popularity, and that its hoped for aim does not seem
so easy to set aside. Namely, as is well known even from common everyday experience,
cooperation will often involve considerable complications and difficulties. First of all, and
already on its most basic {\it conceptual} levels, cooperation proves to be an extremely
complex and rich phenomenon, which therefore cannot in any way be encompassed by a few general
definitions and mathematical models. This fact is indeed in sharp contrast with the modelling
of competition and conflict situations, where for instance in Game Theory, the so called
non-cooperative games, see (2.1) and (2.5) below, describe quite well - and do so in spite of
their manifest simplicity - a considerably large class of such situations. Second, one can
only talk about cooperation if one can {\it rely} on the respective agreements undertaken by
the autonomous agents involved. And the issue of such a reliance clearly depends on a variety
of complex factors which can easily be outside of the realms of convenient mathematical
modelling. \\

Be it as it may, one should not forget that, just like in the case of competition and conflict,
the primary aspect of cooperation is {\it intent}, while the respective subsequent
conceptualizations, models and actions are only specific instances of manifestation and
expression of such an intent. Therefore, the primary issue is whether we intend to have
competition or conflict in a context which may preferably include cooperation as well, or on
the contrary, we intend, because of no matter what reasons, to relegate cooperation to a
secondary role, or even exclude it altogether. \\

And if we do not a priori intend to exclude cooperation, then we should be careful not to
allow that it is excluded merely by default, that is, simply due to the fact that it is in
general not so easy to deal with it, be it conceptually, or practically, and thus it simply
happens that we fail, avoid or decline to make use of it. \\
And a readiness to pursue rational behaviour beyond its non-coopertive limits will then
suggest that the intent to cooperate, and even more importantly, {\it the intent to secure and
keep up in the longer run a context suitable for cooperation} is but a clearly rational
behaviour, even if on a certain meta-level. \\

The effect of the presence of such considerable difficulties related to cooperation has been
that the modelling of cooperation has not received enough attention, see Axelrod. \\

As is well known, Nash himself tended to see Game Theory as being mainly moved by competition,
conflict, and so on, rather than by cooperation, see Nasar. \\

On the other hand, the older and much more experienced John von Neumann, the originator of
modern Game Theory, considered that there was a major and urgent need in economics, and other
important human ventures involving strategic thinking, for the introduction of rational
approaches to the respective variety of human interactions involved. And clearly, the very
attempts to rationalize approaches to competition, conflict, and so on, rather tend to mollify,
than prioritize them, see Neumann \& Morgenstern. \\

As it happened, the first major result in Game Theory was obtained by John von Neumann in his
1928 paper. This is the famous Min-Max theorem about two player zero sum games. Such games
involve the smallest possible nontrivial number of players, and the sharpest possible conflict
among them, in which what one player wins, the other one must lose, thus the sum of what is
won and lost is always zero. Clearly, in such a game there is no way for cooperation. \\

Here it is important to note the following. During the period around 1928, when von Neumann
was only 25 years old, he was involved in at least two other major ventures, namely, the
foundation of Set Theory, and the foundation of Quantum Mechanics, and in both of them he made
most important and lasting contributions. In this way, von Neumann's involvement in Game
Theory during that period can be seen as reflecting the special, if not in fact, fundamental
importance he attributed to it. And indeed, he saw it as being the first ever systematic and
rigorous theoretical approach to a {\it rational} management of conflict and competition
between two or more conscious agents. \\

Problems of optimization had been considered earlier as well. After all, optimization appears
as a rather permanent human concern in a most diverse range of activities. However, such
problems could typically be seen as a game with one single conscious and rational player who
was playing against Nature, or against everybody else. \\
But now, in von Neumann's view, the task was to be able to build an appropriate theory for
conflicts and competitions between a number of conscious agents, assuming that they are firmly
and reliably grounded in rationality. \\
It should be remembered in this regard that von Neumann happened to grow up in the Empire of
Austria-Hungary, and did so during the disastrous years of World War I and its aftermath. And
just like the well known philosopher Karl Popper, of the same generation and background, von
Neumann was much influenced by the prevailing view during the post World War I years, a view
according to which that war - called The Great War - had been the result of nothing else but
systematic and catastrophic, even if rather trivial, failures of rationality on the part of
the leading elites of the Western powers. \\

The extend to which von Neumann gave a special priority to the development of Game Theory is
further illustrated by his activity during the next one and a half decade, till the publishing
in 1944, that is, during World War II, of his joint book with Morgenstern, entitled "Theory of
Games and Economic Behaviour". Indeed, during the years of World War II, von Neumann was
heavily involved in supporting the American war effort and doing so in a variety of ways.
Consequently, at the time, he did very little theoretical research. And yet, he considered it
important enough to dedicate time to Game Theory, and complete the mentioned book of over 600
pages, which is the first ever systematic and detailed presentation of that theory. It should
also be mentioned that the theory in the book is due solely to von Neumann, and most of it,
except for his Min-Max theorem of 1928, was developed by him during the years preceding its
first publication in 1944. Morgenstern was an economist, and his contribution to the book
consisted in the connections between game theory and economic behaviour. In this way, that
book can in fact be seen as a research monograph in Game Theory. \\
As for the relevance of its content, a good deal of it still makes useful reading after more
than six decades. \\

The importance attributed to Game Theory continued after World War II as well. And it was to a
good extent due to the interest manifested in it by the RAND Corporation, a most influential
California think tank at the time, which was heavily involved, among others, in strategic
studies related to the just emerged Cold War. \\

As it happens, Emile Borel initiated in the early 1920s the study of certain well known card
games which were related to the two person zero sum games. However, he did not obtain the
respective major result, and in fact, he assumed that the Min-Max theorem was in general false.
Later, in 1934, R A Fisher was also involved in a study of two person zero sum games, without
however obtaining the major Min-Max theorem, see Luce \& Raiffa. \\

Then starting in 1950, John Nash, who at the time was 22 years old, published his fundamental
papers Nash [1,3] on equilibrium in n-person non-cooperative games, and his main result was a
significant extension of the Min-Max theorem of von Neumann. \\
Nash has also published important results in cooperative games, see Nash [2,4]. \\

However, as it happened, the Nash-equilibrium theorem massively extended the von Neumann
Min-Max theorem, and it does so in two ways. First, it is no longer restricted to two players,
and instead, it can handle an arbitrary number of them. Second, it is no longer restricted to
zero sum games, but it refers to arbitrary non zero sum ones. \\
Furthermore, the Nash-equilibrium theorem has ever since its inception remained the best known
result of its kind, as it has not been further extended in any significant manner, when
considered in its own terms of non-cooperation. \\
Consequently, that theorem has always been seen as paradigmatic for non-cooperative games. \\

To a certain extent, such an interpretation is not so surprising due to the following two
facts. First, the Nash result on the existence of an equilibrium in mixed strategies is an
obvious extension of the corresponding von Neumann Min-Max theorem, and the latter, as
mentioned, is indeed about games which are outside of any possible cooperation. Second, as
long as one is limited to the usual, and thus narrow concepts of cooperation, the result of
Nash on equilibrium will be seen as falling outside of such concepts. \\

As we shall show in section 2, however, such an interpretation can only hold if alone the
usual, and indeed narrow concepts of cooperation are considered. And as the very concept of
equilibrium in the Nash result implies it, that result can have any practical meaning and
value at all, and do so beyond its particular case when only two players are involved, only if
the respective $n \geq$ 3 players do accept - even if implicitly - certain {\it additional
common rules} of behaviour. Thus in the case of $n \geq 3$ players, they {\it must} end up by
cooperating, even if in ways other and more deep than those according to the usual views of
cooperation. And in fact, the kind of cooperation needed in order to enable the
Nash-equilibrium concept and result to function at all proves to be so {\it strong} as to be
practically {\it unrealistic}. \\

Due to the reputation of von Neumann, the interest he showed in Game Theory led in the late
1940s and early 1950s to a considerable status for that theory among young mathematicians at
Princeton, see Nasar. That status was further enhanced by results such as those of Nash and
others. \\
There was also at the time a significant interest in Game Theory outside of academe. As mentioned,
for instance, the well known RAND Corporation was conducting studies in political and military
strategy which were modelled by a variety of games. \\

As it happened, however, soon after, certain major setbacks were experienced. First, and
within game theory itself, was the fact that in the case of n-person games, even for $n \geq$
3 moderately large, there appeared to be serious {\it conceptual} difficulties related to
reasonable concepts of solution. Indeed, too many such games proved not to have solutions in
the sense of various such solution concepts, concepts which seemed to be natural, see Luce \&
Raiffa, Owen, Vorob'ev, Rasmusen. \\
Later, the nature and depth of these conceptual difficulties got significantly clarified. For
instance, in Binmore [1-3] it was shown that there are no Turing machines which could compute
general enough games. In other words, solving games is not an algorithmically feasible
problem. \\

The second major trouble came from outside of Game Theory, and it is not quite clear whether
at the time it was soon enough appreciated by game theorists with respect to its possible
implications about the fundamental difficulties in formulating appropriate concepts of
solutions in games. Namely, Kenneth Arrow showed that a set of individual preferences cannot
in general, and under reasonable conditions, be aggregated into one joint preference, unless
there is a dictator who can impose such a joint preference. This result of Arrow was in fact
extending and deepening the earlier known, so called, voter's paradox, mentioned by the
Marquis de Condorcet, back in 1785, see Mirkin. \\

In subsequent years, developments in Game Theory lost much of their momentum. In later years,
applications of Game Theory gained the interest of economists, and led to a number of new
developments in economic theory. An indication of such developments was the founding in 1989
of the journal Games and Economic Behavior. \\
At the same time certain studies of competition, conflict, and so on, were taken up by the
developments following Arrow's fundamental paper, and led to social, collective or group
choice theory, among others. Decision Theory got also involved in such studies involving
certain specific instances of competition, conflict or cooperation, related to problems of
optimization, see Rosinger [1-6]. \\

In games, or in social, collective or group choice one has many autonomous players,
participants or agents involved, each of them with one single objective, namely, to maximize
his or her advantage which usually is defined by a scalar, real valued utility function. And
by the early 1950s it became clear enough that such a situation would not be easy to handle
rationally, even on a conceptual level. \\
On the other hand, a main objective of Decision Theory is to enable one single decision maker
who happens to have several different, and usually, quite strongly conflicting objectives. And
in view of Arrow's result, such a situation may at first appear to be more easy to deal with.
The situation, however, proves to be quite contrary to such a first perception, see Rosinger
[1-6]. \\

As it turned out, games, social choice or decision making, each have their deeper structural
limitations. \\

And as far as games are concerned, not even in the case of the celebrated and {\it
paradigmatic non-cooperative} Nash-equilibrium theorem is it possible to escape the {\it
paradox} of having the validity of that result essentially based on a {\it cooperative}
assumption so strong as to make it unrealistic in practice. \\ \\

{\bf 2. The Nash-Equilibrium Theorem} \\ \\

The usual way an n-person non-cooperative game in terms of the players' {\it pure strategies}~
is defined is by \\

(2.1) \quad $ G ~=~ ( P,~ ( S_i ~|~ i \in P ),~ ( H_i ~|~ i \in P ) ) $ \\

Here $P$ is the set of $n \geq 2$ players, and for every player $i \in P$, the finite set
$S_i$ is the set of his or her pure strategies, while $H_i : S \longrightarrow \mathbb{R}$ is
the {\it payoff} of that player. Here we denoted by \\

(2.2) \quad $ S ~=~ \prod_{i \in P}~ S_i $ \\

the set of all possible {\it aggregate pure strategies} ~$s = ( s_i ~|~ i \in P) \in S$
generated by the independent and simultaneous individual strategy choices $s_i$ of the players
$i \in P$. \\

The game proceeds as follows. Each player $i \in P$ can freely choose an individual strategy
$s_i \in S_i$, thus leading to an aggregate strategy $s = ( s_i ~|~ i \in P) \in S$. At that
point, each player $i \in P$ receives the payoff $H_i ( s )$, and the game is ended. We assume
that each player tries to maximize his or her payoff. \\ \\

{\bf Remark 1} \\

The usual reason the games in (2.1) are seen as non-cooperative is as follows. Each of the
$n \geq 2$ players $i \in P$ can completely independently of any other player in $P$ choose
any of his or her available strategies $s_i \in S_i$. And the only interaction with other
players happens on the level of payoffs, since the payoff function $H_i$ of the player $i \in
P$ is defined on the set $S$ of aggregated strategies, thus it can depend on the strategy
choices of the other players. \\
However, as we shall see in Remarks 2 - 4 below, in the case of $n \geq 3$ players, this
independence of the players is only apparent, when seen in the framework the concept of
Nash-equilibrium, and the corresponding celebrated Nash theorem.

\hfill $\Box$ \\ \\

Before considering certain concepts of equilibrium, it is useful to introduce some notation.
Given an aggregate strategy $s = ( s_i ~|~ i \in P ) \in S$ and a player $j \in P$, we denote
by $s_{-j}$ what remains from $s$ when we delete $s_j$. In other words $s_{-j} = ( s_i ~|~
i \in P \setminus \{ j \} )$. Given now any $s^{\,\prime}_j \in S_j$, we denote by $( s_{-j},
s^{\,\prime}_j )$ the aggregate strategy $( t _i ~|~ i \in P ) \in S$, where $t_i = s_i$, for
$i \in P \setminus \{ j \}$, and $t_j = s^{\,\prime}_j$, for $i =j$.  \\

For every given player $j \in P$, an obvious concept of {\it best strategy} $s^*_j \in S_j$ is
one which has the {\it equilibrium} property that \\

(2.3) \quad $ H_j ( s_{-j}, s^*_j ) ~\geq~ H_j ( s_{-j}, s_j ),~~~ \mbox{for all}~~
                                                                s \in S,~ s_j \in S_j $ \\

Indeed, it is obvious that any given player $j \in P$ becomes completely independent of all
the other players, if he or she chooses such a best strategy. However, as it turns out, and is
well known, Rasmusen, very few games of interest have such strategies. \\
Consequently, each of the players is in general {\it vulnerable} to the other players, and
therefore must try to figure out the consequences of all the possible actions of all the other
players. \\

Furthermore, even when such strategies exist, it can easily happen that they lead to payoffs
which are {\it significantly lower} than those that may be obtained by suitable cooperation. A
good example in this regard is given by game called the Prisoner's Dilemma, see section
3. \\ \\

{\bf Remark 2} \\

It is precisely due to the mentioned vulnerability of players, which is typically present in
most of the games in (2.1), that there may arise an interest in {\it cooperation} between the
players. A further argument for cooperation comes from the larger payoff individual players
may consequently obtain. A formulation of such a cooperation, however, must then come in {\it
addition} to the simple and general structure present in (2.1), since it is obviously {\it not}
already contained explicitly in that structure.

\hfill $\Box$ \\ \\

Being obliged to give up in practice on the concept of best strategy in (2.3), Nash suggested
the following alternative concept which obviously is much {\it weaker}. \\ \\

{\bf Definition ( Nash )} \\

An aggregate strategy $s^* = ( s^*_i ~|~ i \in P) \in S$ is called a {\it Nash-equilibrium},
if no single player $j \in P$ has the incentive to change {\it all alone} his or her strategy
$s^*_j \in S_j$, in other words, if

\bigskip
(2.4) \quad $ H_j ( s^* ) ~\geq~ H_j ( s^*_{-j}, s_j ),~~~ \mbox{for all}~~ j \in P,~ s_j \in S_j $ \\ \\

{\bf Remark 3} \\

Clearly, the Nash-equilibrium only considers the situation when {\it never more than one
single} player does at any given time deviate from his or her respective strategy. Therefore,
the Nash-equilibrium concept is {\it not} able to deal with the situation when there are $n
\geq$ 3 players, and at some moment, more than one of them deviates from his or her
Nash-equilibrium strategy. \\
Needless to say, this fact renders the concept of Nash-equilibrium {\it unrealistically
particular}, and as such, also {\it unstable} or {\it fragile}. \\

Furthermore, that assumption has a manifestly, even if somewhat implicitly and subtly, {\it
cooperative} nature. \\

Above all, however, the larger the number $n \geq$ 3 of players, the {\it less realistic} is
that assumption in practical cases. \\

It is obvious, on the other hand,  that when there are $n \geq 3$ players, in case at least
two players change their Nash-equilibrium strategies, the game may open up to a large variety
of other possibilities in which some of the players may happen to increase their payoffs. \\

Therefore, when constrained within the context of the Nash-equilibrium concept, the game
becomes {\it cooperative} by {\it necessity}, since the following {\it dichotomy} opens up
inevitably : \\

Either

\begin{itemize}

\item (C1)~~~ All the players agree that never more than one single \\
          \hspace*{1.3cm} player may change his or her Nash-equilibrium strategy,

\end{itemize}

Or

\begin{itemize}

\item (C2)~~~ Two or more players can set up one or more coalitions, \\
          \hspace*{1.3cm} and some of them may change their Nash-equilibrium \\
          \hspace*{1.3cm} strategies in order to increase their payoffs.

\end{itemize}

Consequently, what is usually seen as the essentially non-cooperative nature of the game (2.1),
turns out, when seen within the framework of the Nash-equilibrium concept, to be based - even
if tacitly and implicitly - on the {\it very strong cooperative} assumption (C1) in the above
dichotomy. \\
On the other hand, in case (C1) is rejected, then the game falls out of the Nash-equilibrium
framework, and thus it opens up to the wealth of possibilities under (C2), which among others,
can contain a large variety of possible ways of cooperation. \\

In this way, both the Nash-equilibrium concept and the Nash theorem on the existence of the
respective equilibrium in mixed strategies are highly {\it unstable} or {\it fragile} when
there are 3 or more players involved. \\

Also, similar with the best strategies in (2.3), with the Nash-equilibrium
strategies as well it can happen that they lead to payoffs which are significantly lower than
those that may be obtained by suitable cooperation.

\hfill $\Box$ \\ \\

As in the particular case of (2.1) which gives the von Neumann Min-Max theorem on two person
zero sum games, so with the weakened concept of Nash-equilibrium in (2.4), such an equilibrium
will in general not exist, unless one embeds the {\it pure strategy} game (2.1) into its
extension given by the following {\it mixed strategy} game \\

(2.5) \quad $ \mu G ~=~ ( P,~ ( \mu S_i ~|~ i \in P ),~ ( \mu H_i ~|~ i \in P ) ) $ \\

Here, for $i \in P$, the set $\mu S_i$ has as elements all the probability distributions
$\sigma_i : S_i \longrightarrow [0, 1]$, thus with $\Sigma_{s_i \in S_i}~ \sigma_i ( s_i ) =
1$. Let us now denote \\

(2.6) \quad $ \mu S ~=~ \prod_{i \in P}~ \mu S_i $ \\

Then for $i \in P$ we have the payoff function $\mu H_i : \mu S \longrightarrow \mathbb{R}$
given by \\

(2.7) \quad $ \mu H_i ( \sigma ) ~=~  \Sigma_{s \in S}~ \sigma ( s ) H_i ( s ) $ \\

where for $\sigma = ( \sigma_i ~|~ i \in P ) \in \mu S$ and $s = ( s_i ~|~ i \in P ) \in S$,
we define $\sigma ( s ) = \prod_{i \in P}~ \sigma_i ( s_i )$. \\

Now the definition (2.4) of Nash-equilibrium for pure strategy games (2.1) extends in an
obvious manner to the mixed strategy games (2.5), and then, with the above we have, see
Vorob'ev  \\ \\

{\bf Theorem ( Nash) } \\

The mixed strategy extension $\mu G ~=~ ( P,~ ( \mu S_i ~|~ i \in P ),~ ( \mu H_i ~|~ i
\in P ) )$ of every pure strategy game $G ~=~ ( P,~ ( S_i ~|~ i \in P ),~ ( H_i ~|~ i
\in P ) )$, has at least one Nash-equilibrium strategy. \\ \\

{\bf Remark 4} \\

Obviously, what was mentioned in Remarks 2 and 3 related to the inevitability of cooperation
when the {\it pure strategy} games (2.1) are considered within the framework Nash-equilibrium,
will also hold for the {\it mixed strategy} games (2.5), and thus as well for the above
theorem of Nash. \\ \\

{\bf Remark 5} \\

The idea behind the Nash Program to reduce cooperative games to noncooperative ones seems at
first quite natural. Indeed, in its very essence, a game means that, no matter what the rules
of the game are, each player has a certain {\it freedom} to act within those rules, and can do
so {\it independently} of all the other players. Therefore, it may appear that if we only
concentrate on that freedom and independence, then within that context one can see the game as
noncooperative, that being one of the usual ways to understand the very meaning of freedom and
independence. \\
Furthermore, even if one cooperates, one is still supposed to be left in a game with a certain
freedom and independence. Thus it may still appear that, after subtracting all what is due to
the rules of the game and to one's possible cooperation, one is still supposed to remain with
a certain freedom and independence. \\

According to Nash himself, it could be possible to express all communication and bargaining in
a cooperative game in a formal manner, thus turn the resulting freedom and independence of the
players into moves in an {\it extended} noncooperative game, in which the payoffs are also
extended accordingly. Since such a program has never been fully implemented in all its details
and only its ideas were presented, its criticism must unavoidably remain on the same level of
ideas. However, a certain relevant and well tested objection can be made nevertheless, see
McKinsey [p. 359] : \\

" It is extremely difficult in practice to introduce into the cooperative games the moves
corresponding to negotiations in a way which will reflect all the infinite variety permissible
in the cooperative game, and to do this without giving one player an artificial advantage
( because of his having the first chance to make an offer, let us say )." \\

What is lost, however, in such a view as the Nash Program is that an {\it appropriate
voluntary and mutual limitation} of one's freedom and independence, in order to implement a
{\it cooperation} can {\it significantly change} the payoffs, and thus it can offer to players
an increase in their payoffs, an increase which simply {\it cannot} be attained in any other
noncooperative way. And this is precisely the point in {\it cooperation}. \\

On the other hand, precisely to the extent that the above objection in McKinsey is valid
related to the Nash Program, and all subsequent experience points to its validity, the very
same objection touches essentially on any attempt to reconsider cooperation, and do so in more
formal ways. \\


\begin{thebibliography}{99}

\vspace{0.5cm}

\bibitem{} Arrow, Kenneth J [1] : A difficulty in the concept of social welfare. Journal of
Political Economy, 58, 4 (1950)

\bibitem{} Arrow, Kenneth J [2] : Social Choice and Individual Values, 2nd ed. Wiley, New York, 1963

\bibitem{} Axelrod, R : The Evolution of Cooperation. New York, Basic Books, 1984

\bibitem{} Binmore, Kenneth [1] : Modelling rational players. Part I. Economics and
Philosophy, 3 (1987) 179-214

\bibitem{} Binmore, Kenneth [2] : Modelling rational players. Part II. Economics and
Philosophy, 4 (1988) 9-55

\bibitem{} Binmore, Kenneth [3] : Game theory and the social contract : Mark II
(manuscript 1988) London School of Economics

\bibitem{} Blau, Julian : The existence of social choice functions. Econometrica, 25, 2 (1957) 302-313

\bibitem{} Hargreaves, Shaun P, et. al. : Game Theory, A Critical Intorduction. Routledge,
London, 1995

\bibitem{} Harsanyi, J C : Approaches to the bargaining problem before and after the theory
of games : a critical discussion of Zeuthen's, Hick's, and Nash's theories. Econometrica,
24 (1956), 144-157

\bibitem{} Luce, R Duncan \& Raiffa, Howard : Games and Decisions,
Introduction and Critical Survey. Wiley, New York, 1957, or Dover, New York, 1989

\bibitem{} McKinsey, J C C : Introduction to the Theory of Games. Mc-Graw-Hill, New York,
1952

\bibitem{} Mirkin, Boris G : Group Choice. Wiley, New York, 1979

\bibitem{} Nasar, Silvia : A Beautiful Mind. Faber and Faber, London, 1998

\bibitem{} Nash, John F [1] : Equilibrium points in n-person games. Proc. Nat. Acad. Sci.
USA, 38 (1950), 48-49

\bibitem{} Nash, John F [2] : The bargaining problem. Econometrica, 18 (1950) 155-162

\bibitem{} Nash, John F [3] : Non-cooperative games. Ann. Math., 54 (1951) 286-295

\bibitem{} Nash, John F [4] : Two-person cooperative games. Econometrica, 21 (1953) 128-140

\bibitem{} von Neumann, John : Zur Theorie der Gesellschaftsspiele. Math. Annalen, 100
(1928) 295-320

\bibitem{} von Neuman, John \& Morgenstern, Oskar : Theory of Games and Economic
Behavior. Princeton, 1944

\bibitem{} Owen, Guillermo : Game Theory. Saunders, Philadelphia, 1968

\bibitem{} Rasmusen, Eric : Games and Information. Balckwell, Malden, 2002

\bibitem{} Rosinger, Elemer E [1] : Interactive algorithm for multiobjective optimization.
JOTA, 35, 3 (1981) 339-365

\bibitem{} Rosinger, Elemer E [2] : Errata Corrige : Interactive algorithm for multiobjective
optimization. JOTA, 38, 1 (1982) 147-148

\bibitem{} Rosinger, Elemer E [3] : Aids for decision making with conflicting objectives.
In Serafini, P (Ed.), Mathematics of Multiobjective Optimization. Springer, New York,
1985, 275-315

\bibitem{} Rosinger, Elemer E [4] : Beyond preference information based multiple criteria
decision making. European Journal of Operational Research, 53 (1991) 217-227

\bibitem{} Rosinger, Elemer E [5] : Reconsidering conflict and cooperation.
arXiv:math.OC/0405065

\bibitem{} Rosinger, Elemer E [6] : PIIPTI, or the Principle of Increasing \\ Irrelevance of
Preference Type Information. \\ arXiv:math.OC/0506619

\bibitem{} Tucker, Albert W : A two person dilemma. (unpublished) Stanford University mimeos,
May 1950.

\bibitem{} Vorob'ev, N N : Game Theory, Lectures for Economists and Systems
Scientists. Springer, New York, 1978

\bibitem{} Walker, Paul : An outline of the history of game theory.
(http://william-king.www.drexel.edu/top/class/histf.html)

\bibitem{} Zeuthen, F : Problems of Monopoly and Economic Warfare. Routledge, London,
1930

\end{thebibliography}
\end{document}